\documentclass[11pt,centertags]{amsart}
\usepackage{amssymb}

\usepackage{mathptmx}

\numberwithin{equation}{section}

\copyrightinfo{2006}{V. Kostrykin, K. A. Makarov,  A. K. Motovilov}

\DeclareSymbolFont{SY}{U}{psy}{m}{n}
\DeclareMathSymbol{\emptyset}{\mathord}{SY}{'306}

\newcommand{\bbC}{\mathbb{C}}

\newcommand{\bbR}{\mathbb{R}}
\newcommand{\EE}{\mathsf{E}}

\newcommand{\cK}{{\mathcal K}}

\newcommand{\cO}{{\mathcal O}}

\newcommand{\cU}{{\mathcal U}}

\newcommand{\fH}{\mathfrak{H}}

\newcommand{\fL}{\mathfrak{L}}

\newcommand{\fR}{\mathfrak{R}}

\newcommand{\beq}{\begin{equation}}
\newcommand{\eeq}{\end{equation}}
\newcommand{\ba}{\begin{align}}
\newcommand{\ea}{\end{align}}

\newcommand{\spec}{{\mathrm{spec}}}

\newcommand{\dist}{{\ensuremath{\mathrm{dist}}}}

\newcommand{\diag}{{\ensuremath{\mathrm{diag}}}}

\newcommand{\Ran}{\mathop{\mathrm{Ran}}}

\newtheorem{theorem}{Theorem}[section]
\newtheorem{proposition}[theorem]{Proposition}
\newtheorem{lemma}[theorem]{Lemma}

\newtheorem{hypothesis}[theorem]{Hypothesis}
\theoremstyle{definition}

\newtheorem{example}[theorem]{Example}

\newtheorem{introtheorem}{Theorem}{\bf}{\it}
{\bf}{\it}
{\bf}{\it}
{\bf}{\it}

\theoremstyle{remark}
\newtheorem{remark}[theorem]{Remark}

\begin{document}


\title{Perturbation of spectra and spectral subspaces$^*$}\thanks{$^*$Published
in Trans. Amer. Math. Soc. \textbf{359} (2007), 77 -- 89}

\author[V. Kostrykin]{Vadim Kostrykin}
\address{ Fraunhofer-Institut f\"{u}r
Lasertechnik, Steinbachstra{\ss}e 15, D-52074 Aachen, Germany}
\email{kostrykin@ilt.fraunhofer.de}
\curraddr{Institut f\"ur Mathematik, Technische Universit\"at Clausthal,
Erzstra{\ss}e 1, D-38678 Clausthal-Zellerfeld, Germany}
\email{kostrykin@t-online.de}

\author[K. A. Makarov]{K.~A.~Makarov}
\address{ Department of Mathematics, University of
Missouri, Co\-lum\-bia, MO 65211, USA}
\email{makarov@math.missouri.edu}

\author[A. K. Motovilov]{A.~K.~Motovilov}
\address{ Department of Mathematics, University
of Missouri, Columbia, MO 65211, USA}
\email{motovilv@theor.jinr.ru}
\curraddr{Bogoliubov Laboratory of Theoretical Physics, JINR,
Joliot-Curie 6, 141980 Dubna, Moscow Region, Rus\-sia}

\subjclass[2000]{Primary 47A15, 47A55; Secondary 47B15}


\dedicatory{Dedicated to Volker Enss on the occasion of his 60-th birthday}

\begin{abstract}
We consider the problem of variation of spectral subspaces for linear
self-adjoint operators with emphasis on the case of off-diagonal
perturbations. We prove a number of new optimal results on the shift of the
spectrum and obtain (sharp) estimates on the norm of the difference of two
spectral projections associated with  isolated parts of the spectrum of the
perturbed and unperturbed operators respectively.
\end{abstract}

\maketitle

\section*{Introduction}\label{SMain}
\renewcommand{\theequation}{\arabic{equation}}

It is well known (see, e.g., \cite{Kato}) that if $A$ and
$V$ are bounded self-adjoint operators on a separable
Hilbert space $\fH$, then the spectrum of the operator
$B=A+V$ is confined  in the closed $\|V\|$-neighborhood,
$\cU_{\|V\|}(\spec (A))$, of the spectrum of $A$,
\begin{equation}\label{inclusion1}
\spec(B)\subset \cU_{\|V\|} (\spec (A)).
\end{equation}
In particular,
if the spectrum $\spec(A)$ consists of two isolated parts
$\sigma$ and $\Sigma=\spec(A)\setminus \sigma$ such that
$\dist(\sigma, \Sigma)=d>0$,
under the hypothesis
\begin{equation}\label{natur}
\|V\|<\frac{d}{2}
\end{equation}
the perturbation  $V$ does not close the  gaps in the spectrum of $A$
separating the sets $\sigma$ and $\Sigma$ and hence the spectrum of $B$
also has two separated components. Clearly, condition \eqref{natur} is
sharp in the sense that if  $\|V\|\ge d/2$, the perturbed operator $A+V$
may not have separated parts of the spectrum at all.

One of fundamental problems of the perturbation theory is to study the
variation of the spectral subspace associated with the isolated part
$\sigma$ of the spectrum of $A$ under the perturbation. A still unsolved
problem is to give an answer to the following question: Is it true or not
that under the hypothesis \eqref{natur}
\begin{equation*}
\|\EE_A(\sigma)-\EE_B(\cO_{d/2}(\sigma))\|<1?
\end{equation*}
Here $\EE_T(\Delta)$ denotes the spectral projection for the self-adjoint
operator $T$ corre\-sponding to a Borel set $\Delta$ and $\cO_{d/2}(\sigma)$
is the open $d/2$-neighborhood of the set $\sigma$ (see
\cite{Kostrykin:Makarov:Motovilov:1} for a partially affirmative answer to
this question).

In the present paper we treat the case where the perturbation $V$
is off-diagonal with respect to the direct sum of spectral
subspaces associated with the sets $\sigma$ and $\Sigma=\spec(A)
\setminus \sigma$ in the sense that $$\EE_A(\sigma)
V\EE_A(\sigma)=\EE_A(\Sigma) V\EE_A(\Sigma)=0.$$ \mbox{We address} the
following general question for the off-diagonal perturbations.

\tolerance 500
\begin{itemize}
\item[(i)] What is an optimal  requirement on the norm of
the  perturbation $V$  that \mbox{guarantees} that $V$ does not close the  gaps in
the spectrum of $A$ separating the sets $\sigma$ and $\Sigma$?
\end{itemize}

Unlike the case of general perturbations, in the off-diagonal case the
answer \mbox{depends} upon the mutual disposition of the isolated parts $\sigma$
and $\Sigma$ of the spectrum of the operator $A$. Leaving apart the
well-known case where the sets $\sigma$ and $\Sigma$ are subordinated (see
\cite{Adamyan:Langer:Tretter:2000a}, \cite{Davis:123}, \cite{Davis:Kahan},
\cite{Kostrykin:Makarov:Motovilov:4}) we focus on two cases:

\begin{description}
\item[\textit{Case} $\mathrm{I}$]  {the sets $\sigma$ and $\Sigma$ are separated.\hfill}
\item[\textit{Case} $\mathrm{II}$] {the set $\sigma$ and the convex hull of $\Sigma$ (or vice versa) are
separated.}
\end{description}

We give a complete solution to the problem (i) and show that the
corresponding optimal requirements  are: $\|V\|<\sqrt{3}/2d$ in Case I and
$\|V\|<\sqrt{2} d$ in Case II, respectively.

We also address the following question of perturbation theory for spectral
subspaces.

\begin{itemize}
\item[(ii) ] What can be said about
variation of the spectral subspace associated with the isolated part
$\sigma$ of the spectrum of $A$ under the off-diagonal
perturba\-tions satisfying the optimal requirements above?
\end{itemize}

We \emph{conjecture} that in Case I the inequality $\|V\|<\sqrt{3}/2d$ is
sufficient for the difference of the projections $
\EE_A(\sigma)-\EE_B(\cO_{d/2}(\sigma)) $ to be a strict contraction.  We
also prove that in Case II the optimal ``gap-nonclosing'' requirement
$\|V\|<\sqrt{2} d$ guarantees that
$\|\EE_A(\sigma)-\EE_B(\cO_{d}(\sigma))\|<1$.



\vspace{5mm}

\textbf{Main results.} Let
\begin{equation*}
\delta_V=\|V\| \tan\bigg (\frac{1}{2}
\arctan\frac{2\|V\|}{d}\bigg ).
\end{equation*}
Our \emph{first principal result} is as follows.

\begin{introtheorem}\label{theorem1}
\it Suppose that the self-adjoint bounded perturbation $V$ is
off-diagonal with respect to the decomposition
$\fH=\Ran\EE_A(\sigma)\oplus \Ran\EE_A(\Sigma)$.
If
\begin{equation}\label{normb1}
 \|V\|<\frac{\sqrt{3}}{2}d
\end{equation} or, which is the
same, $\delta_V<d/2$, then

(i) the spectrum of $B$ in the open ${d/2}$-neighborhood $\cO_{d/2}(\sigma)$
of the  set $\sigma$ is separated from the remainder of the spectrum of
$B$. Moreover,
\begin{equation*}
\spec(B) \cap \cO_{d/2}(\sigma)= \spec(B) \cap
\cU_{\delta_V}(\sigma)
\end{equation*}
is a nonempty closed set;

(ii) if in addition $\|V\|< c_\pi d$ with
$c_\pi=\frac{3 \pi - \sqrt{\pi^2 + 32}}{\pi^2 -
4}=0.503288\ldots$, then
\begin{equation*}
\| \EE_A(\sigma)-\EE_{B}(\cO_{d/2}(\sigma))\|\le
\frac{\pi}{2} \frac{\|V\|}{d-\delta_V}<1.
\end{equation*}
\end{introtheorem}

This result is sharp in the sense that if the norm bound
$\|V\|<\frac{\sqrt{3}}{2}d$ is violated, then the set $\spec(B) \cap
\cO_{\delta_V}(\sigma)$  may be either empty or non-closed \, (see
Example \ref{ex:1} below). Theorem \ref{theorem1} implies that the
best possible constant $c$ in the inequality $\|V\|< cd$ implying
$\|\EE_A(\sigma)-\EE_{B}(\cO_{d/2}(\sigma))\|<1$ satisfies the
two-sided estimate
\begin{equation*}
 c_\pi\le c \le\frac{\sqrt{3}}{2},
\end{equation*}
improving the previously known bounds $1/\pi \leq c\leq \sqrt{2}$
\cite{Albeverio} and $\frac{2}{2+\pi}\leq c$
\cite{Kostrykin:Makarov:Motovilov:1}.

If the convex hull  $\cK(\sigma)$ of the set $\sigma$ does not intersect
the remainder $\Sigma$ of the spectrum of $A$ we face a new phenomenon
which does not have an analog in the case of general perturbations. That
is, the spectrum of the component $\Sigma$ may not ``leak out" into the
open $d$-neighborhood   of the set $\sigma$, provided that
$\|V\|<\sqrt{2}d$ and the perturbation $V$ is off-diagonal.

We give a complete solution of the problem in this case
and our \emph{second principal result} is as follows.

\begin{introtheorem}
\label{theorem2} Suppose that the self-adjoint bounded perturbation $V$ is
off-diago\-nal with respect to the decomposition
$\fH=\Ran\EE_A(\sigma)\oplus \Ran\EE_A(\Sigma)$. If $\cK(\sigma)\cap
\Sigma=\emptyset$ and
\begin{equation}\label{normb2}
\|V\|<\sqrt{2}d
\end{equation} or, which is the same, $\delta_V<d$,
then

(i) the spectrum of $B$ in the open {d}-neighborhood
$\cO_d(\sigma)$ of the  set $\sigma$ is separated from the
remainder of the spectrum of $B$. Moreover,
\begin{equation*}
\spec(B) \cap \cO_{d}(\sigma)= \spec(B) \cap
\cU_{\delta_V}(\sigma)\quad \text{ is a nonempty closed
set and}
\end{equation*}

(ii)
\begin{equation*}
\|\EE_A(\sigma)-\EE_{B}(\cO_{d}(\sigma))\|\le \sin\bigg
(\arctan \frac{\|V\|}{d-\delta_V}\bigg ) < 1.
\end{equation*}
\end{introtheorem}

This result is sharp in the following sense. If the norm bound
$\|V\|< \sqrt{2}d $ \mbox{is violated,} then the set $\spec(B) \cap
\cO_{\delta_V}(\sigma)$  may be either empty or non-closed
\mbox{(see \,Example \,\ref{ex:1.6} \,below).} \, Moreover, the best
possible constant $c$ in inequality {$\|V\| < cd$ implying  $\|
\EE_A(\sigma)-\EE_{B}(\cO_{d}(\sigma))\|<1$ is $c=\sqrt{2}$. Note
that the size of the} neighborhood in question is as twice as big as
that in The\-o\-rem~\ref{theorem1}.

As we have already mentioned,  the case where the sets $\sigma$ and
$\Sigma$ are subordinated  is well understood and the following is known
(see \cite{Adamyan:Langer:95}, \cite{Davis:123}, \cite{Davis:Kahan},
\cite{Kostrykin:Makarov:Motovilov:4}).

\begin{introtheorem}\label{2theta} \it
Suppose that the self-adjoint bounded perturbation $V$ is off-diago\-nal
with respect to the decomposition $\fH=\Ran\EE_A(\sigma)\oplus
\Ran\EE_A(\Sigma)$. If the sets $\sigma$ and $\Sigma$ are subordinated and,
for definiteness, $\sup \sigma < \inf \Sigma$, then

(i)  the spectrum of the operator $B$ does not intersect the open interval
$(\sup \sigma,$ $\inf \Sigma)$ and

(ii)
\begin{equation}\label{bla}
\|\EE_{A}(\sigma)-\EE_{B}\bigl((-\infty, \sup
\sigma]\bigr)\|\le \sin \bigg ( \frac{1}{2} \arctan
\frac{2\|V\|}{d}\bigg )< \frac{\sqrt{2}}{2}.
\end{equation}
\end{introtheorem}

In particular, the spectrum of the perturbed operator $B$
always has two subor\-dinated  components and the
perturbation $V$ does not close the gap \, $(\sup \sigma,
\inf \Sigma)$ in the spectrum of $A$ (no requirements on
the norm of $V$ are needed). An analog of Theorem
\ref{2theta} for the case without gap, that is, for $\sup
\sigma \leq \inf \Sigma$  or $\sup \Sigma \leq \inf
\sigma$ is also known (see
\cite{Adamyan:Langer:Tretter:2000a} and
\cite{Kostrykin:Makarov:Motovilov:4}).

A few words about notations. By $\spec (A)$ we denote the
spectrum of a bounded self-adjoint operator $A$ and $\inf
A$ ($\sup A$) denotes the infimum (supremum) of the set
$\spec(A)$. The spectral projection of $A$ associated with
a Borel set $\Delta\subset \bbR$ is denoted by
$\EE_A(\Delta)$ and the resolvent set of $A$ is denoted by
$\rho(A)$. We use the symbol $\cO$ for open sets while the
symbol $\cU$ is usually associated with closed
neighborhoods. If not explicitly stated otherwise, for an
arbitrary orthogonal projection $P$ the symbol $P^\perp$
denotes the orthogonal projection onto the orthogonal
complement of the subspace $\Ran P$ in $\fH$, i.e.,
$P^\perp=I-P$.

\section{Perturbation of Spectra}
\label{Sgeneral}
\numberwithin{equation}{section}

We start this section by presenting a fairly simple but
general result which provides optimal lower and upper
bounds on the shift of the spectrum of a bounded
self-adjoint operator under a perturbation which is
off-diagonal with respect to the given orthogonal
decomposition of the Hilbert space reducing the
unperturbed operator.

\begin{lemma}\label{shura}
Let $A$ and $V$ be  bounded self-adjoint operators on a
Hilbert space $\fH$, $B=A+V$, and $P$ an orthogonal
projection commuting with $A$. Assume, in addition, that
\begin{equation}\label{comm}
P V P=P^\perp V P^\perp=0.
\end{equation}
Denote by $A_0$ and $A_1$ the parts of $A$ associated with
its invariant subspaces $\Ran P$ and $\Ran P^\perp$,
respectively.

Then
\begin{equation}
\label{Bup} \inf A-\delta_V^\ell\le \inf B\le \inf A
\end{equation}
and
\begin{equation}
\label{Bdown} \sup A \leq\sup B \leq \sup A+\delta_V^r,
\end{equation}
where
\begin{align*}
\delta_V^\ell&=\|V\| \tan\left(\frac{1}{2}\arctan
\frac{2\|V\|}
{|\inf A_1-\inf A_0|} \right),\\
\delta_V^r&=\|V\| \tan\left(\frac{1}{2}\arctan
\frac{2\|V\|} {|\sup A_1-\sup A_0|} \right)
\end{align*}
with a natural convention that $ \arctan(+\infty)=\pi/2 $
in the case where $\inf A_1=\inf A_0$ and/or $\sup
A_1=\sup A_0$.
\end{lemma}

\begin{proof} Denote by $W^2(B)$  (cf.\ \cite{LMMT:01}) the
quadratic numerical range of the operator $B$ with respect
to the decomposition $\fH=\Ran P\oplus \Ran P^\perp$,
\begin{equation*}
W^2(B)=\bigcup_{\substack{\|f\|=\|g\|=1 \\ f\in \Ran P\\
g\in \Ran P^\perp}} \spec
\begin{pmatrix}
(f,Bf) & (f,Bg)\\
(g,Bf) & (g,Bg)\\
\end{pmatrix}.
\end{equation*}
For $f$ and $g$ as above taking into account \eqref{comm}
yields
\begin{equation*}
\begin{pmatrix}
(f,Bf) & (f,Bg)\\
(g,Bf) & (g,Bg)\\
\end{pmatrix}=
\begin{pmatrix}
(f,Af) & (f,Vg)\\
(g,Vf) & (g,Ag)\\
\end{pmatrix}=
\begin{pmatrix}
a_0 & v\\
v^* & a_1\\
\end{pmatrix},
\end{equation*}
where we have introduced the notations $a_0=(A_0f,f)$,
$a_1=(g,A_1g)$, and $v=(f,Vg)$. The matrix
$
\begin{pmatrix}
a_0 & v\\
v^* & a_1\\
\end{pmatrix}
$ has two eigenvalues $\lambda$ and $\mu$ given by
\begin{equation*}
\lambda= \min\{a_0, a_1\} -|v|
\tan\left(\frac{1}{2}\arctan \frac{2|v|}{ |a_1- a_0|}
\right)
\end{equation*}
and
\begin{equation*}
\mu=\max\{a_0, a_1\} +|v| \tan\left(\frac{1}{2}\arctan
\frac{2|v|}{ |a_1- a_0|}\right).
\end{equation*}
Clearly the eigenvalues $\lambda$ and $\mu$ satisfy the inequalities
\begin{equation}\label{est1}
\inf A-\delta_V^\ell\le \lambda\le \min\{a_0,a_1\}
\end{equation}
and
\begin{equation}\label{est2}
\max\{a_0,a_1\} \leq\mu \leq \sup A+\delta_V^r.
\end{equation}

Since the quadratic numerical range $W^2(B)$ contains the spectrum
of $B$ while $\inf W^2(B)=\inf B$ and $\sup W^2(B)=\sup B$ (see
\cite{LMMT:01}), estimates \eqref{est1} and \eqref{est2} prove the
assertion, taking into account that $\inf\min\{a_0,a_1\}=\inf A$ and
$\sup\,\max\{a_0,a_1\}=\sup A$.
\end{proof}

Given the result of Lemma \ref{shura}, now the proof of
Theorem \ref{theorem1}\,(i) and Theorem \ref{theorem2}\,(i)
is straightforward.

For notational setup  introduce the following hypothesis.

\begin{hypothesis}\label{h1}
Assume that $A$ and $V$ are bounded self-adjoint operators
on a separable Hilbert space $\fH$. Suppose  that the
spectrum of  $A$  has a part $\sigma$ separated from the
remainder of the spectrum $\Sigma$ in the sense that
\begin{equation*}
 \spec(A)=\sigma\cup\Sigma
\end{equation*}
and
\begin{equation*}
 \dist(\sigma, \Sigma)=d>0.
\end{equation*}
Assume, in addition,  that $V$ is off-diagonal with
respect to the decomposition $\fH=\Ran
\EE_A(\sigma)\oplus\Ran\EE_A(\Sigma)$.
\end{hypothesis}

\begin{theorem}\label{shift}
Assume Hypothesis \ref{h1}. Then

(i) The spectrum of the operator $B$ is contained in the closed
$\delta_V$-neighbor\-hood $\cU_{\delta_V}(\spec (A))$ of the spectrum of $A$
\begin{equation}\label{inclusion}
\spec(B)\subset \cU_{\delta_V} (\spec (A)),
\end{equation}
where
\begin{equation*}
\delta_V= \|V\| \tan\left(\frac{1}{2}\arctan
\frac{2\|V\|}{d}\right ).
\end{equation*}

(ii) If
\begin{equation*}
\|V\|<\frac{\sqrt{3}}{2}d
\end{equation*}
(or, which is  the same, $ \delta_V<d/2$), then the
spectrum of $B$ in the open {d/2}-neighborhood
$\cO_{d/2}(\sigma)$ of the  set $\sigma$ is separated from
the remainder of the spectrum of $B$. That is,
\begin{equation*}
\spec(B) \cap \cO_{d/2}(\sigma)= \spec(B) \cap
\cU_{\delta_V}(\sigma)\quad \text{ is a nonempty closed
set}.
\end{equation*}

(iii) If $\cK(\sigma)\cap \Sigma=\emptyset$ and
\begin{equation}\label{Bq2d}
\|V\|<\sqrt{2}d
\end{equation} (or, which is the same, $\delta_V<d$), then  the
spectrum of $B$ in the open $d$-neighbor\-hood $\cO_d(\sigma)$ of the  set
$\sigma$ is separated from the remainder of the spectrum of $B$. That is,
\begin{equation*}
\spec(B) \cap \cO_{d}(\sigma)= \spec(B) \cap
\cU_{\delta_V}(\sigma)\quad \text{ is a nonempty closed
set}.
\end{equation*}
\end{theorem}

\begin{proof}
(i) Take a $\lambda\in \bbR$  such that
\begin{equation}\label{resol}
\dist\bigl(\lambda, \spec(A) \bigr)>\delta_V.
\end{equation}
Denote by $A_\ell $  the part of the operator $A$
associated with the $A$-invariant subspace
\begin{equation*}
\fL=\Ran \EE_A\bigl((-\infty, \lambda)\bigr)
\end{equation*}
and let $V_{\ell}= \EE_A\bigl((-\infty,
\lambda)\bigr)V|_{\fL}$. By Hypothesis \ref{h1} the
operator  $V$ is off-diagonal with respect to the
decomposition $\fH=\EE_A(\sigma)\fH\oplus\EE_A(\Sigma)\fH,
$ so is $V_\ell$ with respect to the decomposition
$\fL=\EE_{A_\ell}\bigl(\sigma\cap(-\infty,
\lambda)\bigr)\fL\oplus\EE_{A_\ell}\bigl(\Sigma\cap(-\infty,\lambda)\bigr)\fL$.

Applying Lemma \ref{shura} yields
\begin{equation}
\label{avprim} \sup(A_\ell+V_\ell)\leq\sup A_\ell +
\delta_V.
\end{equation}

Similarly introducing $A_r$ as the part of the operator
$A$ associated with the $A$-invariant  subspace
$\fR=\Ran\EE_A\bigl((\lambda,\infty)\bigr)$ and $V_{r}$ as
$\EE_A\bigl((\lambda,\infty)\bigr)V|_{\fR}$ one proves that
\begin{equation}
\label{avpprim} \inf(A_r+V_r)\geq\inf(A_r) - \delta_V.
\end{equation}

Combining \eqref{avprim}, \eqref{avpprim}, and
\eqref{resol} proves that
\begin{equation*}
\sup (A_\ell+V_\ell) <\lambda <\inf(A_r+V_r).
\end{equation*}
Clearly the operator $B$ can be represented as follows
\begin{equation*}
B=\diag\{A_\ell+V_\ell, A_r+V_r\}+W,
\end{equation*}
where $W$ is given by
\begin{equation*}
W=V-\EE_A\bigl((-\infty,\lambda)\bigr)
V\EE_A\bigl((-\infty,\lambda)\bigr) -
\EE_A\bigl((\lambda,\infty)\bigr)
V\EE_A\bigl((\lambda,\infty)\bigr)
\end{equation*}
and $\diag\{A_\ell+V_\ell, A_r+V_r\}$ is a diagonal
$2\times 2$ operator matrix with respect to the
decomposition $\fH=\fL\oplus \fR$. Since $W$ is
off-diagonal with respect to $\fH=\fL\oplus \fR$, and the
spectra of the diagonal entries $A_\ell+V_\ell$ and
$A_r+V_r$ are subordinated,   the
whole interval $(\sup (A_\ell+V_\ell),\inf(A_r+V_r))$
belongs to the resolvent set of $B$ (see, e.g.,
\cite{Davis:Kahan} or \cite{Adamyan:Langer:95}), in particular
$\lambda$ belongs to the resolvent set of the operator
$B$, completing the proof.

Before proving assertions (ii) and (iii) of the theorem,
note that the function
\begin{equation*}
f(x)=x\tan\left(\displaystyle\frac{1}{2}\arctan 2x\right)
\end{equation*}
is  strictly increasing  on the positive semi-axis and, moreover, by direct
computation one gets
\begin{equation*}
f\bigg (\frac{\sqrt{3}}{2}\bigg )=\frac{1}{2}\quad
\text{and}\quad f(\sqrt{2})=1.
\end{equation*}
In particular,
\begin{equation}\label{d/2}
\text{the inequality
}\|V\|<\frac{\sqrt{3}}{2}d\quad\text{implies}\quad
\delta_V<d/2
\end{equation}
and
\begin{equation}\label{d}
\text{the inequality
}\|V\|<\sqrt{2}d\quad\text{implies}\quad \delta_V<d.
\end{equation}

(ii) The part (ii) is an immediate corollary of the part
(i) taking into account \eqref{d/2}.

(iii) Take
\begin{equation*}
\lambda=\sup \sigma +\delta_V
\end{equation*}
and let $A_\ell$, $A_r$, and $V_\ell$, $V_r$ be as above. Note
 that the hypothesis $\cK(\sigma)\cap
\Sigma=\emptyset$ implies $V_r=0$.

Again, as in the proof of (ii)  one concludes that
\begin{equation*}
\sup(A_\ell+V_\ell)\leq\sup\spec(A_\ell) + \delta_V=\sup
\sigma +\delta_V.
\end{equation*}

Hypothesis \eqref{Bq2d} implies that $\delta_V<d$. Since
$V_r=0$, the operator $B$ can be represented in the form
\begin{equation*}
B=\diag\{A_\ell+V_\ell, A_r\}+W,
\end{equation*}
where $W$ is given by
\begin{equation*}
W=V-\EE_A\bigl((-\infty,\lambda)\bigr) V\EE_A\bigl((-\infty,\lambda)\bigr) -
\EE_A\bigl((\lambda,\infty)\bigr) V\EE_A\bigl((\lambda,\infty)\bigr)
\end{equation*}
and $\diag\{A_\ell+V_\ell, A_r\}$ is a diagonal $2\times
2$ operator matrix with respect to the decomposition
$\fH=\fL\oplus \fR$. Since $W$ is off-diagonal with
respect to  $\fH=\fL\oplus \fR$, and the spectra of the
diagonal entries $A_\ell+V_\ell$ and $A_r$ are
subordinated ($\delta_V<d$), the whole interval $(\sup
(A_\ell+V_\ell),\inf(A_r))$ belongs to the resolvent set
of $B$. In particular, the interval $(\sup \sigma
+\delta_V, \sup \sigma +d) $ belongs to the resolvent set
of the operator $B$, that is,
\begin{equation*}
(\sup \sigma +\delta_V, \sup \sigma +d)\subset \rho(B).
\end{equation*}
The proof  of the inclusion
\begin{equation*}
(\inf \sigma-d, \inf \sigma -\delta_V)\subset \rho(B)
\end{equation*}
is analogous.
\end{proof}

\begin{remark}\label{tochArchi}
The  results (ii) and (iii) are optimal. That is, if the
perturbation $V$ is overcritical  in the sense that $\|V\|\ge
\frac{\sqrt{3}}{2}d$ (resp. $\cK(\sigma)\cap \Sigma=\emptyset$ and
$\|V\|\ge \sqrt{2}d$), then the set $\cO_{d/2}(\sigma)\cap \spec
(B)$ (resp. $\cO_{d}(\sigma)\cap \spec (B))$ may be empty.
\end{remark}

The following two  examples illustrate the situation.

\begin{example}\label{ex:1}
Let $\fH=\bbC^4$. Introducing the $4\times 4$  matrices
\begin{equation*}
A=
\begin{pmatrix}
-\frac{3}{2} &  0         &   0 &  0 \\
 0           & -\frac{1}{2}&   0                  & 0\\
 0 &  0  &   \frac{1}{2}       &  0 \\
 0  & 0  &   0                  & \frac{3}{2}
\end{pmatrix}
\quad\text{ and }\quad V=\begin{pmatrix}
0 &       \frac{\sqrt{3}}{2}    &   0 &  0 \\
 0           & 0&   0                  & \frac{\sqrt{3}}{2}\\
 \frac{\sqrt{3}}{2} &  0  & 0      &  0 \\
 0  & 0  &    \frac{\sqrt{3}}{2}
      & 0
\end{pmatrix},
\end{equation*}
one easily verifies that the spectrum of the $4 \times 4$ Jacobi matrix
\begin{equation*}
B=
\begin{pmatrix}
-\frac{3}{2} &      \frac{\sqrt{3}}{2}     &  0 &  0 \\
  \frac{\sqrt{3}}{2}        & -\frac{1}{2}&   0                  & 0\\
 0 &  0  &   \frac{1}{2}       &   \frac{\sqrt{3}}{2}   \\
 0  & 0 &      \frac{\sqrt{3}}{2}              & \frac{3}{2}
\end{pmatrix}
\end{equation*}
consists of the three eigenvalues $-2$, $0$, and $2$, with $0$ being an
eigenvalue of multiplicity two. Setting $\sigma=\{-3/2, 1/2\}$ and
$\Sigma=\{-1/2,3/2\}$ one immediately concludes that in this case
$d=\dist\{\sigma,\Sigma\}=1$ and the perturbation $V$ is off-diagonal with
respect to the decomposition
$\bbC^4=\Ran\EE_A(\sigma)\oplus\Ran\EE_A(\Sigma)$ and $\|V\|=\sqrt{3}d/2$.
However, $\cO_{1/2}(\sigma)=(-2, -1)\cup (0, 1)$ does not intersect the set
$\spec (B)=\{-2, 0, 2\}$.

\end{example}

\begin{example}\label{ex:1.6}
Let $\fH=\bbC^3$,
\begin{equation*}
A =\begin{pmatrix}
-1 &  0         &  0\\
 0           & 0&   0                 \\
0&  0        &  1 \\
 \end{pmatrix}
\quad \text{ and }\quad
 V
=\begin{pmatrix}
0 &   \sqrt{2}      &  0\\
 \sqrt{2}            & 0&   0                 \\
0&  0        &  0 \\
\end{pmatrix}.
\end{equation*}
The spectrum of the $3 \times 3$ matrix
\begin{equation*}
B =\begin{pmatrix}
-1 &    \sqrt{2}        & 0\\
    \sqrt{2}        & 0&   0                 \\
0 &  0        &  1 \\
 \end{pmatrix}
\end{equation*}
consists of the two eigenvalues $-2$ and $1$, with $1$ being an
eigenvalue of multiplicity two. Setting $\sigma=\{0\}$ and
$\Sigma=\{-1,1\}$ one concludes that $\cK(\sigma)\cap
\Sigma=\emptyset$, $d=\dist\{\sigma,\Sigma\}=1$, the perturbation
$V$ is off-diagonal with respect to the decomposition
$\bbC^4=\Ran\EE_A(\sigma)\oplus\Ran\EE_A(\Sigma)$ and
$\|V\|=\sqrt{2}d$. However, $\cO_{1}(\sigma)=(-1, 1)$ does not
intersect the set $\spec (B)=\{-2, 1\}$.

\end{example}

\section{Perturbation of Spectral Subspaces}\label{sec:6}

In this section we accomplish the proof of remaining
statements of Theorem \ref{theorem1} part (ii) and Theorem
\ref{theorem2} part (ii)  related to the perturbation of
spectral subspaces.

\begin{proposition}[\cite{Bhatia:Davis:Koosis},
\cite{Bhatia:Davis:McIntosh}, \cite{McEachin}]\label{mce}
Let $A$ and $B$ be bounded self-adjoint operators and
$\sigma$ and $\Delta$ two Borel sets on the real axis
$\bbR$\,. Then
\begin{equation*}
\dist (\sigma, \Delta) \| \EE_A(\sigma) \EE_B(\Delta)\|
\leq \frac{\pi}{2} \| A-B \|.
\end{equation*}
If, in addition, the convex hull of the set $\sigma$ does
not intersect the set $\Delta$, or the convex hull of the
set $\Delta$ does not intersect the set $\sigma$, then one
has the  stronger result
\begin{equation*}
\dist (\sigma, \Delta) \| \EE_A(\sigma) \EE_B(\Delta)\|
\leq \| A-B \|.
\end{equation*}
\end{proposition}

The proof of Theorem \ref{theorem1} part (ii) is based on
combining Proposition \ref{mce} with information on the
shift of the spectrum  obtained in Theorem \ref{theorem1}
part (i).

\begin{theorem}\label{main}
Assume Hypothesis \ref{h1}. If
\begin{equation}\label{nu}
\|V\|< \frac{3\pi-\sqrt{\pi^2 + 32}}{\pi^2 - 4}d,
\end{equation}
then
\begin{equation*}
\| \EE_A(\sigma)-\EE_{B}(\cO_{d/2}(\sigma))\|\le
\frac{\pi}{2} \frac{\|V\|}{d-\delta_V}<1.
\end{equation*}
\end{theorem}

\begin{proof}
Introduce the notations $P=\EE_{A}(\sigma)$ and
$Q=\EE_{A+V}(\cO_{d/2}(\sigma))$. By Theorem \ref{shift} (i)
\begin{equation*}
Q^\perp=\EE_{B}\big(\cU_{\delta_V}(\Sigma)\big),
\end{equation*}
where $\cU_{\delta_V}(\Sigma)$ denotes the closed
$\delta_V$-neighborhood of the set $\Sigma$.

By the first claim of Proposition \ref{mce},
\begin{equation}\label{pqort}
\|PQ^\perp\|\le \frac{\pi}{2}\frac{\|V\|}{\dist(\sigma,
\cU_{\delta_V}(\Sigma))}.
\end{equation}
The distance between the set $\sigma$ and the
(closed)  $\delta_V$-neighborhood of the set $\Sigma$  can be estimated from
below as follows
\begin{equation*}
\dist(\sigma, \cU_{\delta_V}(\Sigma))\ge d-\delta_V>0
\end{equation*}
using the second claim of Theorem \ref{shift}. Then
\eqref{pqort} implies the inequality
\begin{equation*}
\|PQ^\perp\|\le \frac{\pi}{2}\frac{\|V\|}{d-\delta_V}.
\end{equation*}

It is an elementary exercise to check that the function
\begin{equation*}
f(x)=\frac{\pi}{2}x+x\tan\left(\frac{1}{2}\arctan 2x
\right)-1
\end{equation*}
strictly increases on the positive semi-axis and that $f(x)$
 has a
unique positive root
\begin{equation*}
x=\frac{3 \pi - \sqrt{\pi^2 +
32}}{\pi^2 - 4}.
\end{equation*}
 As a corollary,  hypothesis \eqref{nu} implies the inequality
\begin{equation}\label{gut}
\frac{\pi}{2}\frac{\|V\|}{d-\delta_V}<1.
\end{equation}
Hence,
\begin{equation}\label{pq<1}
\|PQ^\perp\|\le \frac{\pi}{2}\frac{\|V\|}{d-\delta_V}<1.
\end{equation}
Interchanging the roles of $\sigma$ and $\Sigma$ one
obtains the analogous inequality
\begin{equation}\label{qp<1}
\|P^\perp Q\| \le \frac{\pi}{2}\frac{\|V\|}{d-\delta_V}<
1.
\end{equation}
Since
\begin{equation}\label{glaz}
\|P-Q\|=\max \{\|PQ^\perp\|, \|P^\perp Q\|\}
\end{equation}
(see, e.g., \cite[Ch. III, Section 39]{Akhiezer:Glazman}), inequalities
\eqref{pq<1} and \eqref{qp<1} prove the assertion.
\end{proof}

We split the proof of Theorem \ref{theorem2} part (ii)
into several steps.

\textbf{1.} First,  we prove that the difference of the corresponding
spectral projections is a strict contraction.

\begin{lemma}\label{grub}
Assume Hypothesis \ref{h1}. If $\cK(\sigma)\cap
\Sigma=\emptyset$ and
\begin{equation}\label{nunu}
\|V\|<\sqrt{2}d,
\end{equation}
then
\begin{equation}\label{rough}
\| \EE_A(\sigma)-\EE_{B}(\cO_{d}(\sigma))||<1.
\end{equation}
\end{lemma}

\begin{proof}
Introduce the notations $P=\EE_{A}(\sigma)$ and
$Q=\EE_{B}(\cO_{d}(\sigma))$. We also need to introduce
four spectral projections associated with the operators
$A$ and $B$ and let
\begin{equation*}
P_\ell =\EE_{A}\bigl((-\infty, \inf \sigma-d]\bigr)
\quad\text{and}\quad P_r = \EE_{A}\bigl([\sup \sigma+d,
\infty)\bigr)
\end{equation*}
and
\begin{equation*}
Q_\ell=\EE_{B}\bigl((-\infty, \inf \sigma-d]\bigr)
\quad\text{and}\quad Q_r =\EE_{B}\bigl([\sup \sigma+d,
\infty)\bigr).
\end{equation*}

Our first claim is that
\begin{equation*}
\|P_k-Q_k\|<\frac{\sqrt{2}}{2}, \quad k=\ell, r.
\end{equation*}
It can be seen as follows. Since the perturbation $V$ is
off-diagonal with respect to the decomposition
$\fH=\Ran\EE_A(\sigma)\oplus \Ran\EE_A(\Sigma)$,
 the operator $A+V$ can be represented as the
following $3\times 3$ Jacobi type operator-matrix with
respect to the decomposition $\fH=\Ran P_\ell \oplus \Ran
P \oplus \Ran P_r$
\begin{equation*}
B =
\begin{pmatrix}
A_\ell & V_{\ell\sigma} & 0\\
V_{\sigma\ell} & A_{\sigma} & V_{\sigma r} \\
0 & V_{r\sigma}  & A_{r}
\end{pmatrix}.
\end{equation*}
Here we used the notation
\begin{equation*}
A_k=A\vert_{\Ran P_k}, \quad k=\ell, r, \quad
A_\sigma=A\vert_{\Ran P}
\end{equation*}
and
\begin{equation*}
V_{\sigma k}=PV\vert_{\Ran P_k}\quad \text{and}\quad
V_{k\sigma}=V^*_{k\sigma}, \quad k=\ell, r.
\end{equation*}
The perturbation problem $A\longrightarrow B$ can
naturally be split into two subproblems
\begin{equation*}
A=\begin{pmatrix}
A_{\ell} & 0 & 0\\
0 & A_{\sigma} & 0 \\
0 & 0 & A_{r}
\end{pmatrix}
\longrightarrow \widetilde A=
\begin{pmatrix}
A_{\ell} & 0 & 0\\
0 & A_{\sigma} & V_{\sigma r} \\
0 & V_{r\sigma}  & A_{r}
\end{pmatrix}
\longrightarrow B=
\begin{pmatrix}
A_{\ell} & V_{\ell\sigma} & 0\\
V_{\sigma\ell} & A_{\sigma} & V_{\sigma r} \\
0 & V_{r\sigma}  & A_{r}
\end{pmatrix}.
\end{equation*}
The operator matrix $\widetilde A$ is block-diagonal with respect to
the decomposition $\fH=\Ran P_\ell \oplus \Ran P_\ell^\perp$ and
clearly $\|A-\widetilde A\|< \sqrt{2}d$. Applying Theorem
\ref{shift} (i, ii) to the ``lower-dimensional" off-diagonal
perturbation problem
\begin{equation*}
\begin{pmatrix}
 A_{\sigma} & 0 \\
 0 & A_{r}
\end{pmatrix}
\longrightarrow
\begin{pmatrix}
A_{\sigma} & V_{\sigma r} \\
 V_{r\sigma}  & A_{r}
\end{pmatrix}
\end{equation*}
under hypothesis $\|V\|<\sqrt{2}d$ one concludes that the
spectrum of $\widetilde A$ consists of two subordinated
components, $\widetilde \sigma=\spec (A_\ell)=\Sigma\cap
(-\infty, \inf \sigma-d]$ and ``the remainder''
$\widetilde \Sigma$. Moreover,
\begin{equation}\label{subord}
\sup A_\ell= \sup \widetilde \sigma < \inf \sigma
-\delta_V\leq \inf \widetilde \Sigma= \inf\begin{pmatrix}
A_{\sigma} & V_{\sigma r} \\
 V_{r\sigma}  & A_{r}
\end{pmatrix},
\end{equation}
 where
\begin{equation*}
\delta_V=\|V\| \tan\bigg (\frac{1}{2}
\arctan\frac{2\|V\|}{d}\bigg )<d.
\end{equation*}

Applying Theorem \ref{2theta}  to the off-diagonal
perturbation problem $\widetilde A \longrightarrow B$
where the spectra of the diagonal entries $A_\ell$ and
$\begin{pmatrix}
A_{\sigma} & V_{\sigma r} \\
 V_{r\sigma}  & A_{r}
\end{pmatrix}$ are subordinated (cf.\ \eqref{subord}) yields
\begin{equation}\label{pq1}
\|P_\ell-Q_\ell\|<\frac{\sqrt{2}}{2}.
\end{equation}
Using analogous arguments one proves the remaining
estimate
\begin{equation}\label{pq2}
\|P_r-Q_r\|<\frac{\sqrt{2}}{2}.
\end{equation}

Clearly,
\begin{equation*}
\|P^\perp Q\|=\|(P_{\ell}+P_{r})Q\| \leq
\sqrt{\|P_{\ell}Q\|^2+\|P_{r}Q\|^2}
\end{equation*}
and  moreover
\begin{align*}
\|P_{\ell}Q\| \leq \|P_{\ell}(Q+Q_{r})\| &
=\|P_{\ell}Q_{\ell}^{\perp}\|\le\|P_{\ell}-Q_{\ell}\|,\\
\|P_{r}Q\|\leq\|P_{r}(Q+Q_{\ell})\| &
=\|P_{r}Q_{r}^{\perp}\|\le\|P_{r}-Q_{r}\|.
\end{align*}
Thus,
\begin{equation}\label{end1}
\|P^\perp Q\|\leq \sqrt
{\|P_{\ell}-Q_{\ell}\|^2+\|P_{r}-Q_{r}\|^2}<1
\end{equation}
using \eqref{pq1} and \eqref{pq2}. In an analogous way one
proves that
\begin{equation*}
\|PQ^\perp\|<1,
\end{equation*}
and hence
\begin{equation*}
\|P-Q\|=\max\{\|P^\perp Q\|, \|PQ^\perp\|\}<1.
\end{equation*}
The proof is complete.
\end{proof}


\textbf{2.} Next, we obtain the following general result which is of an
\emph{a posteriori} character.
\begin{theorem}\label{tantheta} Assume Hypothesis \ref{h1}.
Suppose, in addition, that $\cK(\sigma)\cap \Sigma=\emptyset$. Denote by
$\cO$ the maximal open (finite or semi-infinite) interval such that $\cO
\cap \Sigma=\emptyset$ and the spectrum $\widetilde\sigma$ of the operator
$B$ in $\cO$,
\begin{equation*}
\widetilde\sigma=\spec (B)\cap \cO,
\end{equation*}
is a closed set.

If
\begin{equation*}
  \|\EE_A(\sigma)-\EE_B(\cO)\|<1,
\end{equation*}
then
\begin{equation}\label{davis}
 \|\EE_A(\sigma)-\EE_B(\cO)\|\le \sin \arctan
 \bigg (\frac{\|V\|}{\dist(\widetilde
 \sigma, \Sigma)}\bigg ).
\end{equation}
\end{theorem}

\begin{proof}
Introduce the notations $P=\EE_A(\sigma)$ and $Q=\EE_B(\cO)$. It is well
known (see \cite[Corollary 3.4]{Kostrykin:Makarov:Motovilov:2}; cf.
\cite[Lemma 2.3]{Apostol:Foias:Salinas}, \cite[Theorem 1]{Daughtry},
\cite{Halmos:69})  that if $\|P-Q\|<1$, then $\Ran Q$ is the graph of a
bounded operator $X: \Ran P \to \Ran P^\perp$ and
\begin{equation}\label{sin}
\|P-Q\|=\frac{\|X\|}{\sqrt{1+\|X\|^2}}.
\end{equation}
Without loss of generality one may assume that
\begin{equation*}
\inf \widetilde \sigma=-\sup \widetilde \sigma.
\end{equation*}
For any $g\in \Ran\EE_B(\cO)$  one obtains that $P^\perp
g=XPg$ and hence
\begin{equation}\label{raz}
P^\perp Bg=P^\perp BPg+P^\perp BP^\perp P^\perp g =P^\perp
BPg + P^\perp B P^\perp X Pg.
\end{equation}
Clearly the following estimates hold
\begin{equation}\label{dva}
\|P^\perp B P g\|=\|V Pg\|\le \|V\|\, \|Pg\|
\end{equation}
and
\begin{equation}\label{tri}
\big (\sup \widetilde \sigma +\dist(\widetilde \sigma,
\Sigma)\big )\|XPg\|\le\|P^\perp B P^\perp X Pg\| .
\end{equation}
If in addition $g\in \Ran Q$ one obtains
\begin{align}\label{cheture}
\|P^\perp Bg\| & =\|P^\perp Q B Q g\|\le \|P^\perp Q\|\,
\|QBQ\|\,\|g\|\nonumber\\
& \le \frac{\|X\|}{\sqrt{1+\|X\|^2}} \sup \widetilde
\sigma \sqrt{1+\|X\|^2}\|Pg\|= \sup \widetilde \sigma
\|X\| \, \|Pg\|.
\end{align}

Combining \eqref{raz}--\eqref{tri}, and
\eqref{cheture} one gets  the inequality
\begin{equation}\label{base}
\bigl(\sup \widetilde \sigma +\dist(\widetilde \sigma,
\Sigma)\bigr)\|XPg\|\le \bigl(\sup \widetilde \sigma
\|X\|+\|V\|\bigr)\|Pg\|, \quad g\in \Ran Q.
\end{equation}
Since \eqref{base} holds for any $g\in \Ran Q$, one concludes that
\begin{equation*}
\|X\|\le \frac{\|V\|}{\dist (\widetilde \sigma , \Sigma)},
\end{equation*}
which proves the assertion in view of \eqref{sin}.
\end{proof}

\begin{remark}
Assertion \eqref{davis} is equivalent to the estimate
\begin{equation*}
\|\tan \Theta\|\le \frac{\|V\|}{\dist (\widetilde \sigma , \Sigma)},
\end{equation*}
 where $\Theta$ is the operator angle between the subspaces
$\Ran\EE_{A}(\sigma)$ and $\Ran\EE_{B}(\cO)$.
(For discussion of this notion see, e.g.,
\cite{Kostrykin:Makarov:Motovilov:2}).
Thus, Theorem \ref{tantheta} is
a generalization of the Davis-Kahan $\tan\Theta$-Theorem which is  one from
four fundamental estimates on the norm of the difference of spectral
projections known as $\sin \Theta$, $ \sin 2\Theta$, $\tan \Theta$, and
$\tan 2\Theta$ Theorems proved by Davis and Kahan in \cite{Davis:123} and
\cite{Davis:Kahan}.
\end{remark}


\textbf{3.} Finally, rough estimate \eqref{rough} of Lemma \ref{grub} can be
sharpened by \emph{a posteriori} result of Theorem \ref{tantheta} in
combination with the result of Theorem \ref{theorem2} part (i). The proof
of Theorem \ref{theorem2} part (ii) is  as follows.

\begin{theorem}
Assume Hypothesis \ref{h1}. If $\cK(\sigma)\cap
\Sigma=\emptyset$ and
\begin{equation}\label{nunu1}
\|V\|<\sqrt{2}d,
\end{equation}
then
\begin{equation}\label{rough1}
\| \EE_A(\sigma)-\EE_{B}(\cO_{d}(\sigma))||\le \sin
\arctan \bigg (\frac{\|V\|}{d-\delta_V}\bigg )<1,
\end{equation}
where
\begin{equation*}
\delta_V=\|V\| \tan\bigg (\frac{1}{2}
\arctan\frac{2\|V\|}{d}\bigg ).
\end{equation*}
\end{theorem}

\begin{proof}
By Theorem \ref{theorem2} part (i) and Lemma \ref{grub}
the set $\cO=\cO_d(\sigma)$ satisfies the hypothesis of
Theorem \ref{tantheta} with
\begin{equation*}
\dist(\widetilde\sigma,\Sigma)=d-\delta_V.
\end{equation*}
Applying Theorem \ref{tantheta} completes the proof.
\end{proof}

\section*{Acknowledgments}

 A.~K.~Motovilov
acknowledges the kind hospitality and financial support by the Department
of Mathematics, University of Missouri, Columbia, MO, USA. He was also
supported in part by the Deutsche Forschungsgemeinschaft and the
Russian Foundation for Basic Research.

\bibliographystyle{amsplain}

\begin{thebibliography}{50}

\bibitem{Akhiezer:Glazman} N.~I.~Achiezer and I.~M.~Glasmann, \textit{Theory of
Linear Operators in Hilbert Space}, Dover Publications,
New York, 1993.
MR1255973 (94i:47001)

\bibitem{Adamyan:Langer:95} V.~Adamyan and H.~Langer, \textit{Spectral properties
of a class of rational operator valued functions}, J.
Operator Theory \textbf{33} (1995), 259 -- 277.
MR1255973 (94i:47001)

\bibitem{Adamyan:Langer:Tretter:2000a} V.~Adamyan, H.~Langer, and C.~Tretter,
\textit{Existence and uniqueness of contractive solutions
of some Riccati equations}, J. Funct. Anal. \textbf{179}
(2001), 448 -- 473.
MR1809118 (2001j:34074)

\bibitem{Albeverio} S.~Albeverio, K.~A.~Makarov, and A.~K.~Motovilov,
\textit{Graph subspaces and the spectral shift function}, Canad. J. Math.
\textbf{55} (2003), 449 -- 503; arXiv: math.SP/0105142.
MR1980611 (2004d:47031)

\bibitem{Apostol:Foias:Salinas} C.~Apostol, C.~Foias, and
N.~Salinas, \textit{On stable invariant subspaces},
Integr. Equat. Oper. Theory \textbf{8} (1985), 721 -- 750.
MR0818331 (87c:47005)

\bibitem{Bhatia:Davis:Koosis} R.~Bhatia, C.~Davis, and P.~Koosis,
\textit{An extremal problem in Fourier analysis with
applications to operator theory}, J. Funct. Anal.
\textbf{82} (1989), 138 -- 150.
MR0976316 (91a:42006)

\bibitem{Bhatia:Davis:McIntosh} R.~Bhatia, C.~Davis, and A.~McIntosh,
\textit{Perturbation of spectral subspaces and solution of
linear operator equations}, Linear Algebra Appl.
\textbf{52/53} (1983), 45 -- 67.
MR0709344 (85a:47020)

\bibitem{Daughtry} J.~Daughtry, \textit{Isolated solutions of
quadratic matrix equations}, Linear Algebra Appl.
\textbf{21} (1978), 89 -- 94.
MR0485926 (58:5720)

\bibitem{Davis:123} C.~Davis, \textit{The rotation of
eigenvectors by a perturbation.  I and II}, J. Math. Anal.
Appl. \textbf{6} (1963), 159 -- 173; \textbf{11} (1965),
20 -- 27.
MR0149309 (26:6799), MR0180852 (31:5082)

\bibitem{Davis:Kahan} C.~Davis and W.~M.~Kahan, \textit{The
rotation of eigenvectors by a perturbation. III}, SIAM J.
Numer. Anal. \textbf{7} (1970), 1 -- 46.
MR0264450 (41:9044)

\bibitem{Halmos:69} P.~R.~Halmos, \textit{Two subspaces}, Trans.
Amer. Math. Soc. \textbf{144} (1969), 381--389.
MR0251519 (40:4746)

\bibitem{Kato} T.~Kato, \textit{Perturbation Theory for Linear
Operators}, Springer--Verlag, Berlin, 1995.
MR0203473 (34:3324)

\bibitem{Kostrykin:Makarov:Motovilov:1} V.~Kostrykin, K.~A.~Makarov,
and A.~K.~Motovilov, \textit{On a subspace perturbation problem},
Proc. Amer. Math. Soc. \textbf{131} (2003), 3469 -- 3476; arXiv:
math.SP/0203240.
MR1991758 (2004c:47029)

\bibitem{Kostrykin:Makarov:Motovilov:2} V.~Kostrykin, K.~A.~Makarov,
and A.~K.~Motovilov, \textit{Existence and uniqueness of solutions to the
operator Riccati eqution. A geometric approach}, in Yu.~Karpeshina,
G.~Stolz, R.~Weikard, Y.~Zeng (Eds.), \textit{Advances in Differential
Equations and Mathematical Physics}, Contemporary Mathematics \textbf{327},
Amer. Math. Soc., 2003, p. 181 -- 198; arXiv: math.SP/0207125.
MR1991541 (2004f:47012)

\bibitem{Kostrykin:Makarov:Motovilov:4} V.~Kostrykin, K.~A.~Makarov,
and A.~K.~Motovilov, \textit{A generalization of the $\tan 2\Theta$
theorem}, in J.~A.~Ball, M.~Klaus, J.~W.~Helton, and L.~Rodman
(Eds.), \textit{Current Trends in Operator Theory and Its
Applications}, Operator Theory: Advances and Applications Vol.~149.
Birkh\"{a}user, Basel, 2004, p.~349 -- 372; arXiv: math.SP/0302020.
MR2063758 (2005d:47041)

\bibitem{LMMT:01} H.~Langer, A.~Markus,  V.~Matsaev,  and C.~Tretter,
\textit{A new concept for block operator matrices: The quadratic numerical
range}, Linear Algebra Appl. \textbf{330} (2001),  89 -- 112.
MR1826651 (2002b:47015)

\bibitem{McEachin} R.~McEachin, \textit{Closing the gap in a
subspace perturbation bound}, Linear Algebra Appl. \textbf{180} (1993), 7 --
15.
MR1206407 (94c:47017)
\end{thebibliography}

\end{document}